\documentclass[12pt]{amsart}
\usepackage{amssymb}
\usepackage[pdftex]{color}

\makeatletter
\def\LaTeX{\leavevmode L\raise.42ex
    \hbox{\kern-.3em\size{\sf@size}{0pt}\selectfont A}\kern-.15em\TeX}
\makeatother

\newcommand{\BibTeX}{{\rm B\kern-.05em{\sc i\kern-.025emb}\kern-.08em\TeX}}

\newcommand{\transpose}{\text{T}}
\newcommand{\supp}{\text{supp}}
\newcommand{\Diff}{\mathcal{D}}
\newcommand{\Diffmu}{\Diff_{\mu}}

\newcommand{\vol}{\text{vol}}

\newcommand{\id}{\text{id}}
\DeclareMathOperator{\Tr}{\text{Tr}}

\newcommand{\grad}{\nabla}
\newcommand{\Laplacian}{\Delta}

\DeclareMathOperator{\diver}{div} \DeclareMathOperator{\curl}{curl}

\newtheorem{theorem}{Theorem}[section]
\newtheorem{proposition}[theorem]{Proposition}

\newtheorem{corollary}[theorem]{Corollary}

\theoremstyle{definition}
\newtheorem{defn}[theorem]{Definition}
\newtheorem{exmp}[theorem]{Example}

\theoremstyle{remark}
\newtheorem{rem}[theorem]{Remark}

\makeatletter
\def\theequation{\thesection.\@arabic\c@equation}
\makeatother

\makeatletter
\def\alphenumi{%
  \def\theenumi{\alph{enumi}}%
  \def\p@enumi{\theenumi}%
  \def\labelenumi{(\@alph\c@enumi)}}
\makeatother

\begin{document}

\title[The WKB method for conjugate points]{The WKB method for conjugate points\\ in the volumorphism group}
\author{Stephen C. Preston}
\address{Department of Mathematics, University of Colorado, Boulder,
CO 80309-0395} \email{Stephen.Preston@colorado.edu}

\thanks{Some of this work was completed while the author was a Visiting Assistant Professor at
Stony Brook University. The author is grateful for their
hospitality.}
\date{\today}

\maketitle

\section{Introduction}\label{introsection}

In this paper, we are interested in the location of conjugate points
along a geodesic in the volumorphism group $\Diffmu(M^3)$ of a
compact three-dimensional manifold $M^3$, without boundary. As shown
in \cite{p}, these are typically pathological, i.e., they can occur
in clusters along a geodesic, unlike on finite-dimensional
Riemannian manifolds. We give an explicit algorithm for finding them
in terms of a certain ordinary differential equation, derived via
the WKB-approximation methods of Lifschitz-Hameiri~\cite{lh} and
Friedlander-Vishik~\cite{fv}. We prove that for a typical geodesic
in $\Diffmu(M^3)$, there will be pathological conjugate point
locations filling up closed intervals in $\mathbb{R}^+$; hence
typically the zeroes of Jacobi fields on the volumorphism group are
dense in intervals.

Interest in the geometry of the volumorphism group dates back to
1966, when Arnold~\cite{a} showed that if the volumorphism group is
equipped with an $L^2$ Riemannian metric, then the geodesics of this
metric correspond to Lagrangian solutions of the ideal fluid
equations. In particular, the right-translated tangent vector to the
geodesic will be an incompressible velocity field satisfying the
Euler equation of ideal fluid mechanics. Arnold then computed the
Riemannian sectional curvature of this metric in some special cases,
showing that it was often negative, and concluded that most fluid
flows were likely to be unstable in the Lagrangian sense. In
addition, he found some sections of positive curvature, and
suggested that the existence of conjugate points would imply that
there were some stable perturbations of a fluid flow.

This geometric approach to equations of continuum mechanics has
since been used by many authors to find partial differential
equations whose solutions represent geodesics on an
infinite-dimensional Lie group. Well-known equations arising in this
way include the Korteweg-de Vries equation, the Camassa-Holm
equation, and the Landau-Lifschitz equation. See for example
\cite{kh} for an overview. In some cases, this approach only works
formally: the geodesic equation resulting may not be a genuine
ordinary differential equation on a Banach manifold in any
reasonable topology. On the volumorphism group, however, the
geodesic equation is very well-behaved in a Sobolev topology: Ebin
and Marsden~\cite{ema} proved in 1970 that if we extend the group of
$C^{\infty}$ volumorphisms $\Diffmu(M)$ to the space of Sobolev
class $H^s$ volumorphisms $\Diffmu^s(M)$, there is a $C^{\infty}$
Riemannian exponential map defined in some open subset of each
tangent space $T_{\id}\Diffmu^s(M)$. In addition, the differential
of the exponential map can be computed in terms of the Jacobi
equation in the same way as on a finite-dimensional manifold;
Misio{\l}ek~\cite{m1} proved in 1992 that the parallel transport and
curvature operators are smooth in the Sobolev topology, so that the
Jacobi equation is also an ordinary differential equation. These
results allow one to study ideal fluids using standard Picard
iteration techniques rather than PDE estimates.

The first examples of conjugate points on the diffeomorphism group
were found by Misio{\l}ek, who in \cite{m1} gave a simple example of
a fluid flow on $\Diffmu(S^3)$ having a conjugate point, and in
\cite{m2} gave a more complicated example of a fluid flow on
$\Diffmu(\mathbb{T}^2)$ having a conjugate point. Arnold and Khesin,
in their comprehensive guide to the geometry of the diffeomorphism
group~\cite{ak}, asked whether one could expect conjugate points to
accumulate along a geodesic, or whether they were isolated as on a
finite-dimensional Riemannian manifold. (In general, as shown by
Grossmann~\cite{g}, conjugate points on an infinite-dimensional
manifold may cluster together along a geodesic. Grossmann defined
monoconjugate and epiconjugate points in terms of the failure of the
differential of the exponential map to be injective or surjective,
respectively, and showed that the two types need not coincide in
infinite-dimensions. In addition, there may be clustering of
monoconjugate or epiconjugate points along geodesics.)

This question was answered by Ebin, Misio{\l}ek, and the
author~\cite{emp}. We showed that for any two-dimensional manifold
without boundary, the exponential map is a nonlinear Fredholm map of
index zero. Thus epiconjugate points and monoconjugate points are
the same, and furthermore all monoconjugate points have finite order
and are isolated along a geodesic. On the other hand, we showed that
on the three-dimensional manifold $D^2 \times S^1$ (the solid torus
with flat metric), a certain rotational flow has an epiconjugate
point that is not monoconjugate. The author~\cite{p} extended this
result to show that on any three-dimensional manifold, the first
conjugate point along any geodesic is pathological: either
epiconjugate and not monoconjugate, or monoconjugate of infinite
order. That paper used the index form to locate the first conjugate
point, and the method was based on finding a simple ordinary
differential equation at each point of the manifold whose solutions
would give an upper bound for the location of a conjugate point.
However, since it relied on the index form rather than the actual
Jacobi equation, it could not be used to find the location of any
conjugate point other than the first.

In this paper, we approximate the Jacobi equation itself by an
ordinary differential equation at each point of the manifold, by
working with highly-peaked perturbations of the fluid flow. This
technique is well-known as the WKB method, first applied in quantum
mechanics. It was first applied to the linearized Euler equation by
Bayly~\cite{b} in a special case, then more generally by
Lifschitz-Hameiri~\cite{lh} and Friedlander-Vishik~\cite{fv} in the
early 90s. By a result of the author~\cite{p1}, the Jacobi equation
can be decomposed into the linearized Euler equation and a simple
linearized flow equation; as a result, we can apply the WKB method
to the Jacobi equation as well. We do this in Section
\ref{mainsection} to obtain the ordinary differential equation
\begin{equation}\label{thebigone}
\pi_{\xi_o^{\perp}} \left[ \frac{d}{dt} \Big( \Lambda(t,x_o)
\frac{d\gamma}{dt}\Big) + \omega_o(x_o) \times
\frac{d\gamma}{dt}\right] = 0, \quad \gamma(0)=0, \quad \langle
\gamma(t),\xi_o\rangle \equiv 0.
\end{equation}
Here $\Lambda(t,x_o)=D\eta(t,x_o)^{\text{T}} D\eta(t,x_o)$ is the
metric deformation of the flow $\eta$ at $x_o\in M$, while
$\omega_o(x_o)=\curl{u}_o(x_o)$ is the initial vorticity at $x_o$.
The vector $\xi_o$ is a fixed unit vector in $T_{x_o}M$. We prove
that if for some fixed $x_o\in M$ and $\xi_o\in T_{x_o}M$, there is
a solution $\gamma(t)$ of \eqref{thebigone} satisfying $\gamma(0)=0$
and $\gamma(a)=0$, then there is an epiconjugate point in
$\Diffmu(M^3)$ located at $t=a$.

The freedom to specify any $x_o$ and any $\xi_o\in T_{x_o}M$ implies
the existence of many epiconjugate points. In Section
\ref{appliedsection}, we consider at a fixed $x_o\in M$ the set of
all $a \in \mathbb{R}^+$ and $\xi_o \in S^2$ such that
\eqref{thebigone} has a solution orthogonal to $\xi_o$ and vanishing
at $t=0$ and $t=a$. By analyzing the equation \eqref{thebigone} in
detail, we show that this set is a collection of two-dimensional
surfaces in $\mathbb{R}^+ \times S^2$, which can be expressed as the
graph of a function on some open subset of $S^2$. From this, we
conclude that the epiconjugate point locations form closed
nontrivial intervals in $\mathbb{R}^+$. As we point out with an
example, these intervals frequently extend to infinity.

\section{Background}

Our basic object is a compact Riemannian manifold $M$ of dimension
$n\ge 2$, without boundary, whose Riemannian metric generates a
volume form $\mu$. (For the most part we will assume $M$ is
three-dimensional, although this is not necessary.) It is a
well-known result of Arnold~\cite{a} that formally the geodesic
equation on the volumorphism group $\Diffmu(M)$, consisting of all
$C^{\infty}$ diffeomorphisms $\eta \colon M \to M$ with $\eta^* \mu
= \mu$, where $\mu$ is the Riemannian volume form, is the equation
of ideal fluid mechanics in Lagrangian coordinates, when
$\Diffmu(M)$ is given the $L^2$ Riemannian metric
\begin{equation}\label{l2metric}
\langle u\circ \eta, v\circ \eta\rangle_{L^2} = \int_M \langle
u,v\rangle \circ \eta \, \mu = \int_M \langle u,v\rangle \,
\mu.\end{equation} Here the tangent space at $\eta\in\Diffmu(M)$ is
given by
$$ T_{\eta}\Diffmu(M) = \{ v \circ \eta \,\vert\,
\diver{v}=0\}.$$ The differentials of the left- and
right-translations are given by
\begin{equation}\label{translations}
DL_{\eta}(v) = D\eta(v) \qquad \text{and} \qquad
DR_{\eta}(v)=v\circ\eta.
\end{equation}
By equation \eqref{l2metric}, the $L^2$ metric is right-invariant
under the group action, but not left-invariant unless $\eta$ is an
isometry.

Although the above works formally, to treat ordinary differential
equations we need a Hilbert manifold structure, which we obtain by
enlarging both $\Diffmu(M)$ and $T_{\eta}\Diffmu(M)$ to the spaces
of Sobolev-class $H^s$ diffeomorphisms and $H^s$ vector fields, with
$s>n/2+1$. Ebin-Marsden established that the geodesic equation for a
diffeomorphism $\eta$ is a smooth second-order differential equation
on this manifold, so that it has smooth solutions and in particular
a smooth exponential map $\exp_{\id}\colon \Omega\subset
T_{\id}\Diffmu^s(M) \to \Diffmu^s(M)$, defined by $\exp_{\id}(u_o) =
\eta(1)$, where $\eta$ is the geodesic with $\eta(0)=\id$ and
$\dot{\eta}(0)=u_o$. (By right-invariance of the metric, it is no
loss of generality to assume all geodesics begin at the identity
map, which we will do from now on.) If $M$ is two-dimensional, then
$\exp_{\id}$ is actually defined on all of $T_{\id}\Diffmu^s(M)$ (a
corollary of a famous result due originally to Wolibner~\cite{w}).
If $M$ is three-dimensional, it is still notoriously unknown whether
$\exp_{\id}$ can be defined on the entire tangent space.

By right-invariance, the geodesic equation may be decomposed into
two equations: the flow equation and the Euler equation. The flow
equation relates material variables to spatial variables, and is
given by \begin{equation}\label{flow} \frac{\partial \eta}{\partial
t}(t,x) = u\big(t, \eta(t,x)\big). \end{equation} The Euler equation
for the velocity field $u$ is
\begin{equation}\label{euler} \frac{\partial u}{\partial t} +
P(\nabla_uu) = 0,
\end{equation}
where $P$ is the $L^2$ orthogonal projection of an arbitrary vector
field onto the space of divergence-free vector fields, given
explicitly by
\begin{equation}\label{projection}
P(w) = w - \grad \Laplacian^{-1}\diver{w}.
\end{equation}
In three dimensions, it is often convenient to rewrite equation
\eqref{euler} (by taking the curl of both sides) as
\begin{equation}\label{curleuler}
\frac{\partial}{\partial t} \curl{u} + [u,\curl{u}] = 0.
\end{equation}

The following was proved by the author in \cite{p}. The main point
is that the linearized Euler equation can be written in two ways,
depending on whether one is right-translating or left-translating
the perturbed tangent vector to the identity. For different
purposes, either equation may be preferred.
\begin{proposition}\label{linearizedeulerprop}
The linearization of equation \eqref{euler}, obtained from a family
$u(t,\sigma)$ of solutions by setting $z(t) =
\partial_{\sigma}\big|_{\sigma=0}u(t,\sigma)$, is
\begin{equation}\label{linearizedeuler}
\frac{\partial z}{\partial t} + P(\nabla_zu + \nabla_uz) = 0
\end{equation}
with $P$ given by \eqref{projection} to keep $z$ divergence-free.
Writing $z = \eta_*v$, we obtain the equation
\begin{equation}\label{leftlinearizedeuler}
P\left(\frac{\partial}{\partial t} \left( \Lambda v\right) +
\omega_o \times v\right) = 0,
\end{equation}
again with the consequence that $v$ is divergence-free. Here
$\Lambda(t,x) = D\eta(t,x)^{\transpose} D\eta(t,x)$ is the metric
pullback and $\omega_o(x) = \curl u(0,x)$ is the initial vorticity.
\end{proposition}

In the study of fluid stability, the linearized Euler equation
\eqref{linearizedeuler} is usually used alone to study Eulerian
stability; however we can also study Lagrangian stability by
incorporating the linearization of the flow equation \eqref{flow}.
From the geometric perspective, this corresponds to studying the
Jacobi equation along the geodesic. The fact that we can decouple
the Jacobi equation this way was originally noticed by
Rouchon~\cite{r}. It has been used to derive many properties of
Jacobi fields and conjugate points; see for example \cite{p1}.

\begin{proposition}\label{jacobisplit}
If $\eta(t)=\exp_{\id}(tu_o)$ is a geodesic in $\Diffmu(M)$ with
velocity field $u(t) = \frac{d\eta}{dt} \circ \eta(t)^{-1}$, then
the Jacobi equation along $\eta$, right-translated to the identity,
may be written as the pair of equations \eqref{linearizedeuler} and
the linearized flow equation \begin{equation}\label{linearizedflow}
\frac{\partial y}{\partial t} + \nabla_uy - \nabla_yu = z.
\end{equation}
Here $\diver{y}=0$ follows from the fact that $\diver{z}=0$.

The differential of the exponential map may be written as
$$ d(\exp_{\id})_{tu_o}(tz_o) = y(t)\circ\eta(t),$$
where $y(t)$ is the solution of equations \eqref{linearizedeuler}
and \eqref{linearizedflow} with initial conditions $y(0)=0$ and
$z(0)=z_o$.
\end{proposition}

\begin{corollary}\label{leftjacobisplitprop}
We may also write the differential of the exponential map as
$$ d(\exp_{\id})_{tu_o}(v_o) = D\eta(t) w(t),$$
where $w(t)$ solves the differential equation
\begin{equation}\label{leftjacobisplit}
P\left(\frac{\partial}{\partial t} \left( \Lambda(t,x)
\frac{\partial w}{\partial t} \right) + \omega_o(x) \times
\frac{\partial w}{\partial t}\right) = 0,
\end{equation}
with initial conditions $w(0)=0$ and $w'(0)=v_o$.
\end{corollary}

\begin{proof}
We have $y(t) = \eta(t)_* w(t)$, and the equation $\frac{\partial
y}{\partial t} + [u,y] = z$ is equivalent (by the definition of the
Lie derivative) to $\frac{\partial}{\partial t}\big(\eta(t)^{-1}_*
y(t)\big) = \eta(t)^{-1}_* z(t)$. Thus we have $\frac{\partial
w}{\partial t} = v$, where $v$ satisfies
\eqref{leftlinearizedeuler}. The corollary follows.
\end{proof}

We know that the exponential map is defined and $C^{\infty}$ on some
neighborhood of the identity in $\Diffmu^s(M)$ for $s>n/2+1$. It is
not defined in the weak $L^2$ topology, since we do not have an
existence and uniqueness theorem for solutions of the Euler equation
with $L^2$ initial velocity. Thus the differential of the
exponential map will only be defined only on some $\Omega\subset
T_{\id}\Diffmu^s(M)\to T_{\exp(tu_o)}\Diffmu^s(M)$. However, if
$u_o$ is sufficiently smooth, we can uniquely extend the
differential to a continuous map on $T_{\id}\Diffmu^0(M)$ in the
$L^2$ topology. (We emphasize that despite the notation,
$\Diffmu^0(M)$ is not a topological group or a smooth manifold in
any known sense.)

\begin{defn}\label{extendeddiff}
Let us denote by $T_{\eta}\Diffmu^0(M)$ the closure of
$T_{\eta}\Diffmu(M)$ in the $L^2$ topology. For any fixed $a\in
\mathbb{R}^+$ and $u_o\in T_{\id}\Diffmu(M)$, we set
$E(a)=(d\exp_{\id})_{au_o}$, and we define $\widetilde{E}(a)\colon
T_{\id}\Diffmu^0(M)\to T_{\exp(au_o)}\Diffmu^0(M)$ to be the
extension of $E(a)$ to the $L^2$-closure of $T_{\id}\Diffmu(M)$.
(This is a closed subspace of the space of all $L^2$ vector fields,
by \cite{ema}.)

We can use $\widetilde{E}$ to extend the notion of conjugate points
to the weak topology defined by the Riemannian metric. We say
$\eta(a)$ is \emph{weakly monoconjugate} to the identity if
$\widetilde{E}(a)$ is not injective, and that $\eta(a)$ is
\emph{weakly epiconjugate} to the identity if $\widetilde{E}(a)$ is
not surjective. Similarly, $\eta(a)$ is \emph{strongly
monoconjugate} if $E(a)$ is not injective, and $\eta(a)$ is
\emph{strongly epiconjugate} if $E(a)$ is not surjective.
\end{defn}

The relations between weakly conjugate points are the same as those
for conjugate points in strong metrics proved by Biliotti, Exel,
Piccione, and Tausk~\cite{bept}, since those authors' proof relies
only on the structure of the Jacobi equation rather than on global
properties of infinite-dimensional manifolds. Thus we have:
\begin{theorem}[Biliotti et al.]\label{weakproposition}
We have the following for weakly conjugate points:
\begin{itemize}
\item The set of weakly monoconjugate points is countable and dense in
the set of weakly epiconjugate points; \item Every weakly
monoconjugate point is also weakly epiconjugate;
\item Any point that is strictly weakly epiconjugate (weakly epiconjugate and not weakly
monoconjugate) is a limit point of weakly monoconjugate points.
\end{itemize}
\end{theorem}

These relations are not in general known to hold for strongly
conjugate points in the present situation, where the Riemannian
$L^2$ metric does not generate the $H^s$ topology in which we are
working. However, we do have the following general results due to
Grossmann~\cite{g}. These results rely only on the symmetry of the
conjugacy relation.

\begin{theorem}[Grossmann]\label{weakstrong}
Weakly and strongly conjugate points both satisfy the following
relationships:
\begin{itemize}
\item Every (weakly/strongly) monoconjugate point is also a
(weakly/strongly) epiconjugate point.
\item If $\eta(a)$ is (weakly/strongly) epiconjugate to the identity, and if
the range of $(\widetilde{E}(a), E(a))$ is closed, then $\eta(a)$ is
also (weakly/strongly) monoconjugate to the identity.
\end{itemize}
\end{theorem}

In addition, the following relationships between weakly and strongly
conjugate points are obvious: every strongly monoconjugate point is
a weakly monoconjugate point (the kernel of the $H^s$ map is a
subset of the kernel of the $L^2$ map), while every strongly
epiconjugate point is a weakly epiconjugate point (the image of the
$H^s$ map is a subset of the image of the $L^2$ map). The following
proposition shows that for smooth geodesics in $\Diffmu(M)$,
strictly weakly epiconjugate implies strongly epiconjugate. For the
basic facts on Fredholm operators that we will use, we refer to
Taylor's Appendix A~\cite{taylor}.

\begin{proposition}\label{weakimpliesstrong}
Suppose the initial velocity $u_o \in T_{\id}\Diffmu^s(M)$ is
$C^{\infty}$, so that the geodesic $\eta(t) = \exp_{\id}(tu_o)$ is
$C^{\infty}$ as long as it exists. If $\eta(a)$ is strictly weakly
epiconjugate to $\id$, then $\eta(a)$ is strongly epiconjugate to
the identity.
\end{proposition}

\begin{proof}
Assume, to get a contradiction, that $\eta(a)$ is strictly weakly
epiconjugate but not strongly epiconjugate to $\id$.

We use the notation of Definition \ref{extendeddiff}. By Proposition
\ref{weakstrong}, $\eta(a)$ is strictly weakly epiconjugate if and
only if the range of $\widetilde{E}(a)$ is not closed, which implies
$\widetilde{E}(a)$ is not Fredholm. In addition, since $\eta(a)$ is
not strongly epiconjugate to $\id$, we know $E(a)$ is surjective,
which means it must also be injective; therefore $E(a)$ is Fredholm.

Thus if we left-translate the Jacobi operators, as
\begin{align*}F(a) &= dL_{\eta(a)^{-1}}\circ E(a) \colon
T_{\id}\Diffmu^s(M) \to T_{\id}\Diffmu^s(M) \\
\widetilde{F}(a) &= dL_{\eta(a)^{-1}}\circ\widetilde{E}(a) \colon
T_{\id}\Diffmu^0(M)\to T_{\id}\Diffmu^0(M),\end{align*} then $F(a)$
is also Fredholm, while $\widetilde{F}(a)$ is not Fredholm. If we
consider in addition the operator
$$\overline{F}(a)=A^s F(a)
A^{-s}\colon T_{\id}\Diffmu^0(M)\to T_{\id}\Diffmu^0(M),$$ where $A
= (1+\curl^2)^{1/2}$, then since $A^s$ is an isomorphism between
$T_{\id}\Diffmu^s(M)$ and $T_{\id}\Diffmu^0(M)$, we know
$\overline{F}(a)$ is also Fredholm. The operators $\widetilde{F}(a)$
and $\overline{F}(a)$ are both defined in the same space, and thus
we can consider their difference.

By Proposition \ref{leftjacobisplitprop}, the operators $F(t)$ and
$\widetilde{F}(t)$ both satisfy the differential equation
\begin{equation}\label{leftequation}
P\left( \frac{d}{dt}\Big( \Lambda(t) \frac{dF}{dt} \Big) + \omega_o
\times \frac{dF}{dt}\right) = 0, \qquad F(0)=0, \quad F'(0) = \id.
\end{equation}
Thus the operator $\overline{F}(t)$ satisfies
$$ P\left(\frac{d}{dt} \Big(\Lambda(t) A^{-s} \frac{d\overline{F}}{dt}
 A^{s} \Big) + \omega_o\Big( A^{-s} \frac{d\overline{F}}{dt}
A^{s}\Big) \right) = 0,$$ where $A = (1+\curl^2)^{1/2}$. By
considering commutators (and using the fact that $[A,P]\equiv 0$),
we have
\begin{equation}\label{commutator}
P\left[ \frac{d}{d t} \Big(\Lambda(t)\frac{d \overline{F}}{d t}\Big)
+ \omega_o \Big( \frac{d \overline{F}}{d t} \Big)\right] = -P\left[
\frac{d}{d t} \Big( [A^s, \Lambda(t)]A^{-s} \frac{d \overline{F}}{d
t} + [A^s, \omega_o] A^{-s} \overline{F}\Big)\right].
\end{equation}

Now since $\Lambda(t)$ and $\omega_o$ are smooth operators (for any
$t$), we know that $[A^s, \Lambda(t)]$ and $[A^s, \omega_o]$ are
both differential operators of order less than $s$, by the product
rule; hence $[A^s, \Lambda(t)]A^{-s}$ and $[A^s, \omega_o]A^{-s}$
are both compact operators on $L^2$. Now given that $\overline{F}$
is a continuous operator from $T_{\id}\Diffmu^0\subset L^2$ to
itself, and given that $P$ is a smooth operator from $L^2$ to
$T_{\id}\Diffmu^0 \subset L^2$, we can write
\begin{equation}\label{hacktacular}
P\left[ \frac{d}{d t} \bigg(\Lambda(t)\Big(\frac{d
(\overline{F}-\widetilde{F})}{d t}\Big)\bigg) + \omega_o
\Big(\frac{d (\overline{F}-\widetilde{F})}{d t}\Big) \right] = K(t)
\end{equation}
where $K(t)$ is a fixed compact operator from $T_{\id}\Diffmu^0
\subset L^2$ to itself. Now $\widetilde{F}(0)=0$ and
$\overline{F}(0)=0$; in addition $\widetilde{F}'(0)=\id$ and
$\overline{F}'(0)=\id$, so that $\widetilde{F}-\overline{F}$
vanishes to both first and second order.

Now consider the equation
$$ P\left[ \frac{d}{dt} \big(\Lambda(t)q(t)\big) + \omega_o q(t)\right] = 0.$$ Since the operators
are all bounded in $L^2$, this is a linear differential equation in
$T_{\id}\Diffmu^0(M)$, with solutions existing for all time; hence
it has a bounded solution operator $G(t,\tau) \colon
T_{\id}\Diffmu^0(M)\to T_{\id}\Diffmu^0(M)$ such that
$q(t)=G(t,\tau)q(\tau)$ for all $t,\tau\in \mathbb{R}$. By Duhamel's
principle, we can then write
$$ \widetilde{F}(t)-\overline{F}(t) = \int_0^t \int_0^{\tau} G(\tau,\sigma)
\circ K(\sigma) \, d\sigma \, d\tau.$$ Since $K(\sigma)$ is compact
for all $\sigma$ and $G(\tau,\sigma)$ is continuous, the composition
is compact. The integrals are limits of sums of compact operators,
hence also compact.

Now by assumption we know that $\overline{F}(t)$ is Fredholm while
$\widetilde{F}(t)$ is not Fredholm; this is a contradiction since
the sum of a Fredholm operator and a compact operator must be
Fredholm.
\end{proof}

We could also try to relate weak monoconjugacy to strong
monoconjugacy, and this could be done most simply if we knew that
every $L^2$ solution to \eqref{leftjacobisplit} with $w(0)=0$ and
$w(T)=0$ were actually $C^{\infty}$ as long as the flow and velocity
field were $C^{\infty}$. (This would imply that any $L^2$
monoconjugate point corresponds to a zero of a $C^{\infty}$ Jacobi
field, and in particular to a zero of an $H^s$ Jacobi field.) While
this is true in many cases, it is not always true, due to the fact
that the order of a monoconjugate point may be infinite.

The example of the $3$-sphere from \cite{p} is illustrative. In this
example, we take a velocity field that is left-invariant under the
standard group action on $S^3$; then it is automatically a steady
solution of the Euler equation and so generates a geodesic of
$\Diffmu(S^3)$. Monoconjugate points occur at $t= \frac{m\pi}{n}$,
where $m$ and $n$ are any positive integers with $m\ge n$. Every
such point has infinite order, and the infinite-dimensional space of
Jacobi fields vanishing at both times is spanned by curl
eigenfields, which are $C^{\infty}$. However, it is easy to find an
infinite sum of such fields which converges in $L^2$ but not in
$H^s$ for any $s>0$, just by choosing the coefficients correctly.
Hence there is a Jacobi field vanishing at $t=0$ and
$t=\frac{m\pi}{n}$ which is in $L^2$ but not in $H^s$ for any $s>0$.

In spite of this, it is still true that every weak monoconjugate
point along this particular geodesic in $\Diffmu(S^3)$ is actually a
strong monoconjugate point, and we believe it is likely that this is
always true. However the example above shows the proof may be
somewhat subtle.

The following theorem from \cite{p} tells us when the first weakly
conjugate point occurs, in terms of an ordinary differential
equation.

\begin{theorem}\label{oldmain}
Suppose $\eta\colon [0,T) \to \Diffmu(M)$ is a geodesic. Let $u$ be
the velocity field defined by $\frac{\partial \eta}{\partial t}(t,x)
= u\big(t,\eta(t,x)\big)$. Let us define $\Lambda(t,x) =
D\eta(t,x)^{\transpose} D\eta(t,x)$ and $\omega_o(x) =
\curl{u_o}(x)$.

For each $x\in M$, let $\tau(x)>0$ be the first time such that
\begin{equation}\label{oldmainequation}
\frac{d}{dt} \Big( \Lambda(t,x) \frac{du}{dt}\Big) + \omega_o(x)
\times \frac{du}{dt} = 0 \end{equation} has a solution vanishing at
$t=0$ and $t=\tau(x)$. Then the first weakly conjugate point to the
identity along $\eta$ occurs at $\inf_{x\in M} \tau(x)$. This point
is either strictly weakly epiconjugate (i.e., weakly epiconjugate
but not weakly monoconjugate), or it is weakly monoconjugate of
infinite order.
\end{theorem}

\section{Locating epiconjugate points}\label{mainsection}

Since the Jacobi equation is closely related to the linearized Euler
equation by Proposition \ref{jacobisplit}, we can use techniques
developed by Lifschitz and Hameiri~\cite{lh} and Friedlander and
Vishik~\cite{fv} for the latter. These authors showed that solutions
of the linearized Euler equation for sharply peaked initial data
could be approximated near the peak by a certain ordinary
differential equation, using a WKB approximation. Lifschitz-Hameiri
proved a weaker estimate that is valid more generally than that of
Friedlander-Vishik, but one that is sufficient for our purposes.
(Friedlander-Vishik's estimate is valid only for steady solutions of
the Euler equation on the flat torus $\mathbb{T}^3$, while
Lifschitz-Hameiri is valid for all solutions in any geometry.)

The basic technique of WKB analysis is to posit a solution in the
form $$ z = e^{i\Phi/\varepsilon} a + r,$$ with $\varepsilon$ a
small parameter, and expand in powers of $\varepsilon$ to obtain
simple equations for $\Phi$ and $a$. Then we prove that the error
term $r$ is bounded by $\varepsilon$. The following computation and
error estimate was performed by Lifschitz and Hameiri~\cite{lh}.

\begin{theorem}[Lifschitz-Hameiri]\label{wkbapproximation}
Suppose $u$ is a (possibly time-dependent) solution of the Euler
equation \eqref{euler} on a compact manifold $M$, on some time
interval $[0,T]$.

If $z = e^{i\Phi/\varepsilon} a + r$ solves the linearized Euler
equation \begin{equation}\label{linearizedeulerexplicit}
\frac{\partial z}{\partial t} + \nabla_uz + \nabla_zu = 2 \grad
\Laplacian^{-1}\diver{(\nabla_zu)},\qquad \diver{z} = 0,
\end{equation} with initial condition $z(0) = e^{i\Phi_o/\varepsilon} a_o$,
then for $\varepsilon$ small, the dominant terms satisfy
\begin{equation}\label{Phiequation} \frac{\partial \Phi}{\partial t}
+ u(\Phi) = 0,
\end{equation} and \begin{equation}\label{Aequation} \frac{\partial
a}{\partial t} + \nabla_ua + \nabla_au = 2q \grad \Phi,
\end{equation}
where \begin{equation}\label{requation} q = \frac{\langle \nabla_au,
\grad \Phi\rangle}{\langle \grad \Phi, \grad
\Phi\rangle},\end{equation} while the remainder term $r$ satisfies
\begin{equation}\label{l2error}
\lVert r(t)\rVert_{L^2} \le C\varepsilon
\end{equation}
for all $t\in [0,T]$, where $C$ is a constant that depends on $T$
and derivatives of $\Phi$ and $a$.
\end{theorem}

Equations \eqref{Phiequation} and \eqref{Aequation} look like
partial differential equations, but they are more properly thought
of as ordinary differential equations, which makes analyzing their
solutions simpler. We can write them using right-translation as ODEs
along a particular Lagrangian path, or alternatively using
left-translation as ODEs in a particular tangent space.

\begin{proposition}\label{odes}
If $u$ is a solution of the Euler equation \eqref{euler} with flow
$\eta$, then we can write the solution of equation
\eqref{Phiequation} as $\Phi(t) = \Phi_o\circ\eta(t)^{-1}$ with
$\grad \Phi(t) = (D\eta(t)^{-1})^{\star} \grad \Phi_o$. In addition,
we have for each $x$ the formula $a\big(t,\eta(t,x)\big) =
\alpha(t)$ where $\alpha(t)$ is a vector field along $t\mapsto
\eta(t,x)$ satisfying the equation
\begin{equation}\label{oderight}
\frac{D\alpha}{dt} + \nabla_{\alpha(t)} u(t) = 2\,
\frac{\langle\nabla_{\alpha(t)}u(t), \grad \Phi(t)\rangle}{\langle
\grad \Phi(t), \grad \Phi(t)\rangle} \grad \Phi(t), \quad
\alpha(0)=a(0,x), \quad \langle \alpha(t), \grad
\Phi(t)\rangle\equiv 0.
\end{equation}

Alternatively, for each $x_o$ we can write $a\big(t,\eta(t,x)\big) =
D\eta(t,x)\big(\beta(t)\big)$, where $\beta(t)\in T_{x_o}M$ for all
$t$ and $\beta(t)$ satisfies
\begin{equation}\label{odeleft}
\frac{d}{dt} \big(\pi_{\xi_o^{\perp}} \Lambda(t) \beta(t)\big) +
\langle \xi_o, \omega_o\rangle \xi_o \times \beta(t) = 0, \quad
\beta(0)= a(0,x),\quad  \langle \beta(t), \xi_o\rangle\equiv 0,
\end{equation}
where $\xi_o = \grad \Phi_o(x)$ and $\omega_o = \curl{u_o}(x)$.
\end{proposition}

\begin{proof}
The equation $\frac{\partial \Phi}{\partial t} + u(\Phi)=0$ can be
rewritten as $\frac{\partial}{\partial t} \big(
\Phi(t)\!\circ\!\eta(t)\big) = 0$, which implies the solution
$\Phi(t) = \Phi_o\circ\eta(t)^{-1}$. The formula for $\grad \Phi(t)$
is a consequence of the chain rule. Finally the formula
$$ \frac{D}{dt}a\big(t,\eta(t,x)\big) = \frac{\partial a}{\partial t}\big(t,\eta(t,x)\big) +
\nabla_{u(t,\eta(t,x))}a\big(t,\eta(t,x)\big) $$ allows us to write
\eqref{Aequation} as the ordinary differential equation
\eqref{oderight}.

To obtain \eqref{odeleft} from \eqref{oderight}, we apply
$(D\eta)^{\transpose}$ to both sides of \eqref{oderight} and obtain
\begin{equation}\label{preprojection}
(D\eta)^{\transpose}\left(\frac{D\alpha}{dt} +
\nabla_{\alpha}u\right) = c \xi_o,
\end{equation}
for some function $c$.
The left side was computed in \cite{p} to be, with $\alpha =
D\eta(\beta)$,
$$ (D\eta)^{\transpose}\left(\frac{D\alpha}{dt} + \nabla_{\alpha}u\right) =
\frac{d}{dt} \big( \Lambda(t) \beta(t)\big) + \omega_o \times
\beta(t),$$ using conservation of vorticity. Applying
$\pi_{\xi_o^{\perp}}$ to both sides of \eqref{preprojection} gives
\eqref{odeleft}.
\end{proof}

\begin{corollary}\label{leftl2estimate}
If $w(t)$ is a solution of equation \eqref{leftjacobisplit} with
initial condition $w(0,x)\equiv 0$ and $w_t(0,x) =
e^{i\Phi_o(x)/\varepsilon} v_o(x)$, then $w(t,x) =
e^{i\Phi_o(x)/\varepsilon} \int_0^t \beta(\tau,x)\,d\tau  +
\tilde{r}(t,x)$, where $\beta(t,x)$ satisfies \eqref{odeleft} with
initial condition $\beta(0,x)=v_o(x)$ and $\lVert
\tilde{r}(t)\rVert_{L^2} \le \tilde{C}\varepsilon$ for $t\in [0,T]$.
Here $\tilde{C}$ a constant depending on $T$ and derivatives of
$\eta$ and $v_o$. \end{corollary}

\begin{proof}
We have $w(t,x) = \int_0^t v(\tau, x)\,d\tau$ and $$v(t) =
\eta(t)^{-1}_* z(t) = e^{i\Phi_o/\varepsilon} \eta(t)^{-1}_* a(t) +
\eta(t)^{-1}_* r(t),$$ where $a$ satisfies \eqref{Aequation}. Thus
\begin{align*}\lVert w(t) - e^{i\Phi_o/\varepsilon}
\beta(t)\rVert_{L^2} &\le \int_0^t \lVert \eta(\tau)^{-1}_*
r(\tau)\rVert \, d\tau \\
&\le C\varepsilon \int_0^t \lVert \eta(\tau)^{-1}_* \rVert_{L^2}
\,d\tau \\
&\le C\varepsilon \int_0^t \lVert \eta(\tau)^{-1} \rVert_{C^1} \,
d\tau
\end{align*}
since $\eta$ is volume-preserving. With $\tilde{C} = C\int_0^t\lVert
\eta(\tau)^{-1} \rVert_{C^1} \, d\tau$, we are done.
\end{proof}

The main theorem of the paper follows.

\begin{theorem}\label{main}
Let $\eta$ be a smooth geodesic in $\Diffmu(M)$ with Eulerian
velocity field $u$. If for some $x_o\in M$ and some unit-length
$\xi_o\in T_{x_o}M$, the equation
\begin{equation}\label{itsnotreallyselfadjoint}
\frac{d}{dt}\left( \pi_{\xi_o^{\perp}} \Lambda(t, x_o) \,
\frac{d\gamma}{dt}\right) + \langle \omega_o(x_o), \xi_o\rangle
\xi_o \times \frac{d\gamma}{dt} = 0, \qquad \langle \gamma(t),
\xi_o\rangle\equiv 0\end{equation} has a solution with $\gamma(0)=0$
and $\gamma(a)=0$, then $\eta(a)$ is weakly epiconjugate to the
identity along $\eta$.
\end{theorem}

\begin{proof}
Let $\delta>0$ be any number small enough that Riemannian normal
coordinates exist in a $\delta$-neighborhood of $x_o$. Choose normal
coordinates $(x_1,x_2,x_3)$ with axes aligned so that $\xi_o =
\partial_{x_1}\big|_{x_o}$ and $\gamma'(0)=
\partial_{x_3}\big|_{x_o}$, where $\gamma$ is the solution of
\eqref{itsnotreallyselfadjoint} with $\gamma(0)=\gamma(a)=0$. Define
$\Phi_o \colon B_{\delta}(x_o) \to \mathbb{R}$ by
$\Phi_o(x_1,x_2,x_3)=x_1$, and set $\xi_o \equiv \grad \Phi_o =
\partial_{x_1} + O(\delta^2)$. (We could extend $\Phi_o$ globally
with a bump function, but this is not necessary, since we will only
be working in this neighborhood.) Let $\psi \colon \mathbb{R} \to
\mathbb{R}$ be a bump function with $\psi\equiv 1$ in
$[-\frac{1}{2}, \frac{1}{2}]$, $\psi\equiv 0$ on $\mathbb{R}
\setminus (-1,1)$, and $\lvert \psi'\rvert \le 3$ everywhere. Define
a function $\alpha \colon M \to \mathbb{R}$ by $\alpha(x_1,x_2,x_3)
= x_2 \psi\left( \frac{x_1^2+x_3^2}{\delta^2} +
\frac{x_2^2}{\delta^4}\right)$ (and zero outside the coordinate
neighborhood). Finally, define $v_o = \xi \times \grad \alpha$. Then
in the support of $v_o$, we must have $x_1<\delta$, $x_2<\delta^2$,
and $x_3<\delta$. As a result, we have
\begin{multline*} \grad \alpha = \frac{2x_1x_2}{\delta^2} \, \psi'\left(
\frac{x_1^2+x_3^2}{\delta^2} +
\frac{x_2^2}{\delta^4}\right)\,\partial_{x_1} +
\frac{2x_3x_2}{\delta^2} \, \psi'\left( \frac{x_1^2+x_3^2}{\delta^2}
+ \frac{x_2^2}{\delta^4}\right)\,\partial_{x_3} \\ +
\left[\psi\left( \frac{x_1^2+x_3^2}{\delta^2} +
\frac{x_2^2}{\delta^4}\right) + \frac{2x_2^2}{\delta^4}\psi'\left(
\frac{x_1^2+x_3^2}{\delta^2} +
\frac{x_2^2}{\delta^4}\right)\right]\,\partial_{x_2} +
O(\delta^2),\end{multline*} and since $x_1x_2=O(\delta^3)$ and
$x_3x_2=O(\delta^3)$, we have
$$ \grad \alpha = \left[\psi\left( \frac{x_1^2+x_3^2}{\delta^2} +
\frac{x_2^2}{\delta^4}\right) + \frac{2x_2^2}{\delta^4}\psi'\left(
\frac{x_1^2+x_3^2}{\delta^2} +
\frac{x_2^2}{\delta^4}\right)\right]\,\partial_{x_2} + O(\delta).$$
Therefore, $$v_o(x_1,x_2,x_3)=\left[ 2\frac{x_2^2}{\delta^4}
\psi'\left( \frac{x_1^2+x_3^2}{\delta^2} +
\frac{x_2^2}{\delta^4}\right) + \psi\left(
\frac{x_1^2+x_3^2}{\delta^2} + \frac{x_2^2}{\delta^4}\right)\right]
\,
\partial_{x_3} + O(\delta).$$
Since $x_2 = O(\delta^2)$, $\psi = O(1)$, and $\psi' = O(1)$, both
terms in $v_o$ are $O(1)$ in the support of $v_o$.

Let us write $\Gamma(t,x) \colon \xi_o(x)^{\perp} \subset T_xM \to
\xi_o(x)^{\perp}$ for the solution operator of
$$ \frac{\partial}{\partial t} \left( \pi_{\xi_o(x)^{\perp}} \Lambda(t,x) \,
\frac{\partial \gamma}{\partial t}(t,x) \right) + \frac{\langle
\omega_o(x), \xi_o(x)\rangle}{\langle \xi_o(x), \xi_o(x)\rangle} \,
\xi_o(x) \times \frac{\partial \gamma}{\partial t}(t,x) = 0$$ with
initial conditions $\Gamma(0,x) = 0$ and $\partial_t \Gamma(0,x) =
\id$, i.e., $\Gamma(t,x)(\chi_o)=\gamma(t)$ where $\gamma(0)=0$ and
$\dot{\gamma}(0)=\chi_o$. Then $\Gamma(t,x)$ is smooth and does not
depend on $\delta$; thus we can write $\Gamma(a, x) = \Gamma(a, x_o)
+ O(\delta)$. By assumption and our coordinate construction, we have
$\Gamma(a, x_o)(\partial_z) = 0$. Since $\Gamma(t,x_o)$ is an
operator only in the space perpendicular to $\xi_o=\partial_{x_1}
\vert_{x_o}$, and since $\Gamma(a,x_o)(\partial_{x_3}) = 0$ by
construction of our coordinates, we have for any $\chi\in
\xi_o(x_o)^{\perp}$ that $\Gamma(a, x_o)(\chi) = \langle \chi,
\partial_{x_2}\rangle \Gamma(a,
x_0)(\partial_{x_2})+ O(\delta)$.

In particular for the $v_o$ constructed above, we have
\begin{equation*} \Gamma(a, x)\big(v_o(x)\big) = \langle v_o(x),
\partial_{x_2}\rangle \Gamma(a,x_o)(\partial_{x_2})  + O(\delta) =
O(\delta), \end{equation*} so that $$ \lVert
\Gamma(a,x)\big(v_o(x)\big)\rVert_{L^2} = O(\delta) \cdot \vol
[\supp (v_o)] = O(\delta^5). $$ On the other hand, we have $ \lVert
v_o \rVert_{L^2} = O(\delta^4)$ since $v_o = O(1)$ on $\supp(v_o)$.

By Corollary \ref{leftl2estimate}, we have
$$ \lVert w(a)\rVert_{L^2} \le B \delta^5 + \tilde{C} \varepsilon$$
for some constant $B$. Since $\tilde{C}$ depends on derivatives of
$v_o$, it may have a $\delta$-dependence; however we are still free
to choose $\varepsilon$, and thus we can choose it small enough so
that $\lVert w(a)\rVert_{L^2} = O(\delta^5)$. On the other hand, we
still have $\lVert w'(0)\rVert_{L^2} = O(\delta^4)$.

Now we consider a normalized sequence $w_n$ of solutions with
$\delta_n = \frac{1}{n^2}$ and $\varepsilon_n$ chosen so that
$\lVert w_n(a)\rVert_{L^2} = O(\frac{1}{n^2})$ and $\lVert
w_n'(0)\rVert_{L^2} = 1$. Then the series $\sum_{n=1}^{\infty}
w_n(a)$ converges to some divergence-free vector field in $L^2$. We
now have two possibilities: either $\widetilde{E}(a)$ is injective,
or it is not. If $\widetilde{E}(a)$ is not injective, then $\eta(a)$
is weakly monoconjugate to $\id$ and thus also weakly epiconjugate
to $\id$, and we are done. Otherwise, if $\widetilde{E}(a)$ is
injective, then we have $\sum_{n=1}^{\infty}
\widetilde{E}(a)\big(w_n'(0)\big)$ convergent in $L^2$, while
$\sum_{n=1}^{\infty} w_n'(0)$ cannot converge in $L^2$, so that
$\widetilde{E}(a)$ does not have closed range. Hence $\eta(a)$ must
be weakly epiconjugate to $\id$. Thus we are done.
\end{proof}

\begin{rem}
In the proof of Theorem \ref{main}, we assumed for simplicity that
the dimension of $M$ is $3$. However the theorem holds generally as
long as $\dim{M}\ge 2$. In the more general case we need to replace
the operator $\chi \mapsto \curl{u_o} \times \chi$ with $\chi
\mapsto (\iota_{\chi} du_o^{\flat})^{\sharp}$, where the generalized
vorticity cross product is defined in arbitrary dimensions to
satisfy (for any $\zeta$)
$$\langle (\iota_{\chi} du_o^{\flat})^{\sharp}, \zeta \rangle =
dX_o^{\flat}(\chi,\zeta) = \langle \nabla_{\chi}u_o,\zeta\rangle -
\langle \nabla_{\zeta}u_o,\chi\rangle.$$ With this modification, all
other parts of the proof are valid, with the expression
$x_1^2+x_3^2$ replaced whenever it appears with $x_1^2$ (in
dimension two) or $x_1^2 + x_3^2+ x_4^2 + \cdots + x_n^2$ (in
dimension $n\ge 3$).

However, in dimension two the theorem happens not to give any useful
information. This is because we must have $\langle
\gamma,\xi_o\rangle \equiv 0$, which in two dimensions implies that
$\gamma = f \xi_o^{\perp}$ for some function $f$. (Here
$\xi_o^{\perp}$ is the standard rotation by $90^{\circ}$.) As a
result, $dX_o^{\flat}(\gamma,\dot{\gamma}) = f\dot{f}
dX_o^{\flat}(\xi_o^{\perp}, \xi_o^{\perp}) = 0$ regardless of the
function $f$, so that equation \eqref{itsnotreallyselfadjoint}
becomes
$$ \frac{d}{dt}\left( \langle \xi_o^{\perp}, \Lambda(t, x_o) \xi_o^{\perp}\rangle \,
\frac{df}{dt}\right) = 0. $$ This equation clearly has \emph{no}
solutions $f(t)$ vanishing at two distinct times, due to
positive-definiteness of $\Lambda$, and hence Theorem \ref{main}
does not yield any epiconjugate points at all in dimension two.

This is as expected, since the exponential map is Fredholm in
dimension two~\cite{emp}. If there were nontrivial epiconjugate
points obtainable by Theorem \ref{main}, it would imply (as shown in
the next section) that epiconjugate points could occur in intervals
along a geodesic; but Fredholmness of the exponential map implies
that conjugate points must be isolated along a geodesic.

On the other hand, the fact that the same technique works and yields
genuine epiconjugate points in any dimension higher than three
implies that Fredholmness must also fail in any dimension higher
than three.
\end{rem}

Although Theorem \ref{main} allows us to locate many of the
conjugate points along a particular geodesic, it will typically not
get all of them. Essentially, the reason for this is the following:
the first conjugate point is obtained from the solution of a
three-dimensional second-order differential equation at some point,
by Theorem \ref{oldmain}, while the epiconjugate points detectable
by the WKB method all come from the solution of a two-dimensional
second-order differential equation. In essence, we are looking at
the vanishing of the index form
\begin{equation}\label{pointwiseindex}
I_{0,\tau}(\gamma,\gamma) = \int_0^{\tau} \left\langle \Lambda(t,x)
\frac{d\gamma}{dt}, \frac{d\gamma}{dt}\right\rangle + \left\langle
\omega_o(x)\times \gamma(t), \frac{d\gamma}{dt}\right\rangle \, dt
\end{equation}
among vectors $\gamma(t)$ vanishing at $t=0$ and $t=\tau$. Both
problems have the same index form, but the WKB problem forces the
vectors $\gamma(t)$ to all be orthogonal to the same fixed vector.
Thus there are fewer solutions.

Explicitly, suppose the first conjugate point location $\tau$ along
a geodesic is found from equation \eqref{oldmainequation}, such that
that equation has only one solution vanishing at both $t=0$ and
$t=\tau$, with nontrivial variation in all three directions. Then no
two-dimensional vector solution can make the index form vanish, and
thus the first epiconjugate point coming from Theorem \ref{main}
must be beyond the first conjugate point coming from Theorem
\ref{oldmain}. We construct such an example in the following.

\begin{exmp}
Consider the annular solid torus $S^1 \times S^1 \times [a,b]$ (the
region enclosed by two concentric tori) with coordinate system
$(x,y,z)$ in which the metric is given by
$$ ds^2 = z \, dx^2 + \big(dy + f(x)\,dz\big)^2 + dz^2,$$
for some periodic function $f(x)$. If $u = \partial_x$ in this
coordinate system, then $\curl{u} = \partial_y$, so that we
automatically have $[u,\curl{u}]=0$, which implies by
\eqref{curleuler} that $u$ is a steady solution of the Euler
equation.

In the orthonormal basis $e_1 = z^{-1/2} \, \partial_x$, $e_2 =
\partial_y$, $e_3 = z^{1/2} (\partial_z - f(x)\,\partial_y)$,
the operators in \eqref{oldmainequation} take the form
$$ \Lambda(t,x,y,z) = \left(\begin{matrix} 1 & 0 & 0 \\ 0 & 1 & k(t,x,z)
\\ 0 & k(t,x,z) & 1 + k(t,x,z)^2\end{matrix} \right)\qquad
\text{and}\qquad \omega_o(x,y,z) = \left(\begin{matrix} 0 & 0 & 1 \\
0 & 0 & 0 \\ -1 & 0 & 0 \end{matrix}\right),$$ where $k(t,x,z) =
\sqrt{z} \big(f(t+x)-f(x)\big)$. The solution operator $\Upsilon$ of
\eqref{oldmainequation} satisfying $\Upsilon(0)=0$ and
$\Upsilon'(0)=\id$ is
\begin{equation}\label{oldmainsolnop}
\Upsilon(t) = \left( \begin{matrix} \sin{t}  & F(t)\sin{t}-G(t)\cos{t} & \cos{t}-1 \\
-G(t) & t+\int_0^t k(\tau)^2\,d\tau +\int_0^t \big[G(\tau)F'(\tau)-F(\tau)G'(\tau)\big]\,d\tau & -F(t) \\
1-\cos{t} & -F(t)\cos{t}-G(t)\sin{t} & \sin{t}\end{matrix}\right)
\end{equation}
where $F(t) = \int_0^t k(\tau)\cos{\tau}\,d\tau$ and $G(t)=\int_0^t
k(\tau)\sin{\tau}\,d\tau$ (suppressing the dependence on $x$ and $z$
for simplicity).

There is a conjugate point at $t$ iff there is some vector $v_o$
with $\Upsilon(t)(v_o) = 0$, i.e., iff $\det{\Upsilon(t)}=0$. We can
easily compute
\begin{multline*} \det{\Upsilon(t)} = -\sin{t} \Big[F(t)^2+G(t)^2\Big] \\+
2(1-\cos{t}) \left[ t+\int_0^t k(\tau)^2\,d\tau +\int_0^t
\big[G(\tau)F'(\tau)-F(\tau)G'(\tau)\big]\,d\tau\right].\end{multline*}
For a particular example, e.g., $f(x) = \sin{x}$, it is easy to
check numerically that for any $(x,y,z)$, the first time where
$\det{\Upsilon(t)}=0$ is $t=2\pi$. (This is probably true in
general, but it is not important here).

The unique (up to rescaling) vector $v_o$ for which
$\Upsilon(2\pi)(v_o)=0$ is
$$ v_o = \left( \begin{matrix} F(2\pi) \\ 0 \\
-G(2\pi)\end{matrix}\right),$$ and the corresponding solution to
\eqref{oldmainequation} is
$$ \gamma(t) = \left( \begin{matrix} F(2\pi)\sin{t} + G(2\pi)(1-\cos{t})
\\
G(2\pi)F(t)-F(2\pi)G(t) \\
F(2\pi)(1-\cos{t})-G(2\pi)\sin{t}\end{matrix}\right).$$ There is no
vector $\xi_o$ for which $\langle \gamma(t),\xi_o\rangle = 0$ for
all $t$. As a result, the reduced index form does not vanish at
$t=2\pi$, and so the WKB method does not predict the epiconjugate
point at $t=2\pi$.
\end{exmp}

\section{Conjugate point intervals}\label{appliedsection}

We now explore the consequences of Theorem \ref{main}. First, we
have an easy theorem on continuity of conjugate point locations with
respect to the fixed initial vector $\xi_o$.

\begin{theorem}\label{continuousdependence}
If for some $x_o \in M$ and some unit vector $\xi_o \in T_{x_o}M$,
the equation \eqref{itsnotreallyselfadjoint} has a solution
$\gamma(t)$ with $\gamma(0)=0$ and $\gamma(t_o)=0$, then for any
$\tilde{\xi}_o$ sufficiently close to $\xi_o$, there is a solution
$\tilde{\gamma}(t)$ with $\tilde{\gamma}(0)=0$ and
$\tilde{\gamma}(\tilde{t}_o)=0$ for some $\tilde{t}_o$.
\end{theorem}

\begin{proof}
Fix an orthonormal basis of vector fields $e_1$, $e_2$, and $e_3$ in
a neighborhood of $x_o$. Using spherical coordinates $\theta$ and
$\phi$, we can write for any $x$ near $x_o$ \begin{align*}\xi_o &=
\sin{\theta}\cos{\phi} e_1\big|_{x} + \sin{\theta}\sin{\phi}
e_2\big|_{x} + \cos{\theta}
e_3\big|_{x} \\
\xi_1 &= \cos{\theta}\cos{\phi} e_1\big|_{x} +
\cos{\theta}\sin{\phi}e_2\big|_{x} -\sin{\theta} e_3\big|_{x} \\
\xi_2 &= -\sin{\phi} e_1\big|_{x} + \cos{\phi} e_2\big|_{x}.
\end{align*}
Then any vector orthogonal to $\xi_o$ must be a linear combination
of $\xi_1$ and $\xi_2$, so that the equation
\eqref{itsnotreallyselfadjoint} becomes, with $u(t) = f(t) \xi_1 +
g(t) \xi_2$, \begin{multline}\label{itsmoreselfadjoint}\frac{d}{dt}
\left[ \left(\begin{matrix} \langle \xi_1, \Lambda(t,x) \xi_1\rangle
& \langle \xi_1, \Lambda(t,{x})\xi_2\rangle
\\ \langle \xi_1, \Lambda(t,{x})\xi_2\rangle & \langle \xi_2,
\Lambda(t,{x}) \xi_2\rangle \end{matrix}\right)
\left(\begin{matrix} \dot{f}(t) \\
\dot{g}(t)\end{matrix}\right)\right]\\ + \langle \omega_o(x) \times
\xi_1, \xi_2\rangle \, \left(
\begin{matrix} 0 & 1
\\ -1 & 0 \end{matrix} \right) \left(\begin{matrix}
\dot{f}(t) \\ \dot{g}(t)\end{matrix}\right)  = \left( \begin{matrix} 0 \\
0
\end{matrix} \right).
\end{multline}
In this way, we can consider the dependence of the conjugate time
$t_o$ on $\theta$, $\phi$, and $x$.

By the general theory of oscillation for self-adjoint systems (see
Reid~\cite{reid}), equation \eqref{itsmoreselfadjoint} has a
solution satisfying $u(0)=0$ and $u(t_o)=0$ if and only if the index
form
\begin{multline}\label{indexformpoint}
I_T(f\xi_1 + g\xi_2,f\xi_1+g\xi_2) =\\
\int_0^T \Big( \tilde{\Lambda}_{11}(t,\theta,\phi) \dot{f}(t)^2 +
2\tilde{\Lambda}_{12}(t,\theta,\phi) \dot{f}(t)\dot{g}(t) +
\tilde{\Lambda}_{22}(t,\theta,\phi) \dot{g}(t)^2 +
2\tilde{\omega}_o(\theta,\phi) f(t) \dot{g}(t)\Big) \, dt,
\end{multline}
defined for vector functions $u(t)$ vanishing at $t=0$ and $t=T$, is
negative on some subspace for $T>t_o$ and positive-definite for
$T<t_o$. (Here, for each $x\in M$, we have
$\tilde{\Lambda}_{ij}(t,\theta,\phi) = \langle \xi_i,
\Lambda(t)\xi_j\rangle$ and $\tilde{\omega}_o(\theta,\phi) = \langle
\omega_o \times \xi_1, \xi_2\rangle$.)

In fact, if  we set $$ J(\theta,\phi)(\gamma) =
\frac{I_T(\gamma,\gamma)}{\int_0^T \lvert \gamma(t)\rvert^2 \,
dt},$$ then for any fixed $\theta_o$ and $\phi_o$, we will have for
any $T<t_o$ that the infimum satisfies  $\inf_{\gamma}
J(\theta_o,\phi_o)(\gamma) = \lambda(T,\theta_o,\phi_o)> 0$, while
if $T>t_o$, then $\inf_{\gamma} J(\theta_o,\phi_o)(\gamma) =
\lambda(T,\theta_o,\phi_o)<0$. We can prove this by noting that, by
standard Sturm-Liouville theory for self-adjoint systems,
$J(\theta,\phi)(\gamma)$ attains a minimum at a certain eigenvector
field $\gamma(t)$ of a self-adjoint operator, and the minimum is an
eigenvalue which we can denote by $\lambda(T,\theta,\phi)$. This
eigenvalue depends continuously on $T$, $\theta$, and $\phi$ by the
usual theory.

As a result, for any fixed $T<t_o$, if $\theta$ and $\phi$ are
sufficiently close to $\theta_o$ and $\phi_o$, we will have
$\lambda(T,\theta,\phi)>0$; and for any fixed $T>t_o$ we will have
$\lambda(T,\theta,\phi)<0$. As a result, by continuity with respect
to $T$ we must have $\lambda(\tilde{t}_o,\theta,\phi)=0$ for some
$\tilde{t}_o$.
\end{proof}

The theorem above gives existence of (possibly multivalued)
functions $t$ of spherical variables. For any fixed $\theta_o$ and
$\phi_o$ we have an open domain in the sphere, containing $\theta_o$
and $\phi_o$, on which $t(\theta,\phi)$ is defined. Hence
$t(\theta,\phi)$ is defined on some open subset of the $2$-sphere.
We can understand its behavior better by examining the differential
equation directly, as follows. As a consequence, we will prove that
$t$ is a differentiable function of $\theta$ and $\phi$.

\begin{proposition}\label{traceprop}
At a point $x_o\in M$ and for a unit $\xi_o \in T_{x_o}M$ with
$\langle \xi_o, \omega_o(x_o)\rangle \ne 0$, the equation
\eqref{itsnotreallyselfadjoint} has a solution with $\gamma(0)=0$
and $\gamma(a)=0$ if and only if the solution of the $2\times 2$
matrix equation
\begin{equation}\label{matrixeq} \frac{dW}{dt} + \langle \omega_o(x_o),
\xi_o\rangle J \big(\pi_{\xi_o^{\perp}}
\Lambda(t,x_o)\pi_{\xi_o^{\perp}}\big)^{-1} W = 0
\end{equation}
with $W(0)=\id$ satisfies $\Tr{W(a)} = 2$. Here $J =
\left(\begin{smallmatrix} 0 & -1 \\ 1 & 0 \end{smallmatrix}\right)$.
\end{proposition}

\begin{proof}
Let us write the solution $\gamma(t)$ of
\eqref{itsnotreallyselfadjoint} with $\gamma(0)=0$ and
$\gamma'(0)=v_o$ as a matrix operator $u(t) = S(t)(v_o)$, with $S(t)
\colon \xi_o^{\perp} \to \xi_o^{\perp}$. Then the $2\times 2$ matrix
$\Gamma$ must satisfy
\begin{equation}\label{2ndmatrixeq}
\frac{d}{dt} \left(\pi_{\xi_o^{\perp}} \Lambda \,
\frac{dS}{dt}\right) + \langle \omega_o, \xi_o\rangle J\,
\frac{dS}{dt} = 0, \end{equation} with $S(0)=0$ and
$\frac{dS}{dt}(0) = \id$.

Writing $W = \pi_{\xi_o^{\perp}} \Lambda \, \frac{dS}{dt}$ and using
the fact that $\Lambda(0)=\id$, we have the initial condition
$W(0)=\id$, while $W$ satisfies \eqref{matrixeq}. In addition, since
$\pi_{\xi_o}^{\perp} \Lambda \pi_{\xi_o}^{\perp}$ is symmetric and
$J$ is antisymmetric, we have $\Tr \frac{dW}{dt} W^{-1} = 0$, so
that $\det{W(t)} \equiv 1$ for all $t$. On the other hand, if we
integrate \eqref{2ndmatrixeq} in time, we obtain
$$ W(t) + \langle \xi_o, \omega_o\rangle J S(t) = \id.$$
Since by assumption $\langle \xi_o, \omega_o\rangle \ne 0$, there is
a $v_o \in \xi_o^{\perp}$ with $S(a)(v_o)=0$ if and only if $W(a)$
has $1$ as an eigenvalue. Since $W$ is a $2\times 2$ matrix with
$\det W = 1$, the eigenvalues of $W$ satisfy $\lambda^2 - \Tr{W(a)}
\lambda + 1 = 0$, so that $\lambda=1$ is a solution if and only if
$\Tr{W(a)}=2$.
\end{proof}

Now to study differentiability, we fix the dependence of the
solution on $\xi_o\in S^2$ as in Theorem \ref{continuousdependence};
then equation \eqref{2ndmatrixeq} becomes
\begin{equation}\label{simplematrix} \frac{\partial
W(t,\theta,\phi)}{\partial t} + J \Theta(t,\theta,\phi)
W(t,\theta,\phi) = 0,
\end{equation}
with $\Theta(t,\theta,\phi)=\langle \omega_o(x_o), \xi_o\rangle
(\pi_{\xi_o^{\perp}} \Lambda(t,x_o)\pi_{\xi_o^{\perp}})^{-1}$. In
this way, we get a function $W(t,\xi_o)$ defined on the $2$-sphere.

\begin{theorem}\label{solvefort}
Let $W(t,\xi_o)$ be the solution of equation \eqref{simplematrix}
with $W(0,\xi_o)=\id$, and suppose we have $\Tr{W(t_o, \xi_o)}=2$
for some $t_o>0$ and $\xi_o\in S^2$ with $\langle \xi_o,
\omega_o(x_o)\rangle\ne 0$. Then we can solve the equation
$\Tr{W(t,\xi_o)} = 2$ for $t$ in terms of $\xi$, differentiably in
any neighborhood of $\xi_o$.
\end{theorem}

\begin{proof}
Let us consider the $2$-sphere as a manifold; it is enough to prove
that if $y$ is a coordinate on $S^2$ in some system, then
$\frac{\partial t}{\partial y}$ exists. So let us suppose the other
coordinate as fixed, and think of $W(t,y)$ as depending only on $y$.
We assume $y=y_o$ corresponds to the specified point $\xi_o$ in
$S^2$.

Since $\det W(t,y)\equiv 1$, we can show using the Cayley-Hamilton
Theorem that we have $\big[ W(t_o,y_o) - \id\big]^2 = 0$. Thus
$W(t_o,y)-\id = N$ for some nilpotent matrix $N$, which must look
like $N = cJ\zeta\zeta^T$ for some unit vector $\zeta$ and real
number $c$. Directly from equation \eqref{simplematrix}, we have
\begin{align*}
\frac{\partial}{\partial t} \Tr W(t_o, y_o) &= -\Tr \big( J
\Theta(t_o,y_o) (\id + c
J\zeta\zeta^T)\big)\\
 &= -c\Tr \big( J\Theta(t_o,y_o)
J\zeta\zeta^T\big) \\&= -c\Big( \langle \zeta, J\Theta(t_o,y_o)
J\zeta\zeta^T (\zeta)\rangle + \langle J\zeta, J\Theta(t_o,
y_o) J\zeta\zeta^T (J\zeta)\rangle\Big)\\
&= c \langle J\zeta, \Theta(t_o,y_o) J\zeta\rangle.
\end{align*}
In the first line, we used the familiar fact that if $A$ is
antisymmetric and $B$ is symmetric, then $\Tr{AB} = 0$. (We will use
this again many times throughout the proof.) Now the term $\langle
J\zeta, \Theta(t_o, y_o) J\zeta\rangle$ is positive since $\Theta$
is positive-definite, so that if $c\ne 0$, then $\frac{\partial \Tr
W(t_o,y_o)}{\partial t} \ne 0$. As a result, we can solve for $t$ in
terms of $y$ by the implicit function theorem (and hence obtain
$\frac{d t}{d y}$) if $W(t_o,y_o) \ne \id$. On the other hand, if
$W(t_o,y_o)=\id$, then $\frac{\partial \Tr W}{\partial t}(t_o,y_o) =
0$. Thus we have to work a bit harder.

So now we suppose that $W(t_o,y_o)=\id$. Since the coefficients of
\eqref{simplematrix} are smooth in $y$, we can differentiate
$W(t,y)$ with respect to $y$. If we differentiate
\eqref{simplematrix} with respect to $y$ and write $\frac{\partial
W}{\partial y} = WM$, then $M$ must satisfy the equation
\begin{equation}\label{Peq}
\frac{\partial M}{\partial t} = - W^{-1} J \frac{\partial
\Theta}{\partial y} W,
\end{equation}
and since $\Tr (W^{-1} J \frac{\partial \Theta}{\partial y} W) = \Tr
(J\frac{\partial\Theta}{\partial y}) \equiv 0$, we have
$\frac{\partial}{\partial t} \Tr M(t,y) \equiv 0$ for all $t$ and
$y$, so that $\Tr M(t,y) = 0$ for all $t$ and $y$. Thus since
$W(t_o,y_o)=\id$, we have $\frac{\partial \Tr W}{\partial
y}(t_o,y_o) = \Tr M(t_o, y_o) = 0$.

So both partial derivatives of $\Tr W$ are zero at $(t_o,y_o)$, and
thus we can  look at second derivatives. If $\frac{dt}{dy}$ exists,
then we must have (by differentiating $\Tr W\big(t(y),y\big)$ twice
with respect to $y$)
$$ \frac{\partial^2 \Tr W}{\partial t^2}(t_o,y_o)\, \left(\frac{dt}{dy}\right)^2 +
2 \frac{\partial^2 \Tr W}{\partial t\partial y}(t_o,y_o) \,
\frac{dt}{dy} +  \frac{\partial^2 \Tr W}{\partial y^2}(t_o,y_o)=0,$$
and in order to be able to solve this for $\frac{dt}{dy}$, we must
have
\begin{equation}\label{sufficientimplicit}
\left(\frac{\partial^2 \Tr W}{\partial t\partial
y}(t_o,y_o)\right)^2 - \frac{\partial^2 \Tr W}{\partial
t^2}(t_o,y_o) \, \frac{\partial^2 \Tr W}{\partial y^2}(t_o,y_o) \ge
0.
\end{equation} Conversely if this condition is satisfied, we obtain
either one or two values for $\frac{dt}{dy}$, corresponding to a
crossing of at most two solutions.

Observe first that, quite generally, if $W$ is a matrix with $\det W
\equiv 1$, then for any parameter $r$ we have $\Tr (W^{-1}
\tfrac{\partial W}{\partial r}) = 0$. If we then differentiate with
respect to any other parameter $s$, we have
$$ \Tr (W^{-1} \tfrac{\partial^2 W}{\partial r\partial s}) = \Tr
(W^{-1} \tfrac{\partial W}{\partial r} W^{-1} \tfrac{\partial
W}{\partial s}).$$ Now in our special case, we have $W(t_o,y_o)=\id$
so that
\begin{equation}\frac{\partial^2 \Tr W}{\partial\theta\partial t}(t_o,y_o)
= \Tr \left(\frac{\partial W}{\partial t}(t_o,y_o) \, \frac{\partial
W}{\partial y}(t_o,y_o)\right)
\end{equation}
while
\begin{equation}\frac{\partial^2 \Tr W}{\partial y^2} (t_o,y_o) =
\Tr \Big( \frac{\partial W}{\partial y}(t_o,y_o)\Big)^2 = -2
\det\Big( \frac{\partial W}{\partial y}(t_o,y_o)\Big)
\end{equation}
and
\begin{equation}\frac{\partial^2 \Tr W}{\partial t^2}(t_o,y_o) =
\Tr \Big( \frac{\partial W}{\partial t}(t_o,y_o)\Big)^2 = -2
\det\Big( \frac{\partial W}{\partial t}(t_o,y_o)\Big)
\end{equation}
(these last two are consequences of the Cayley-Hamilton theorem).

Thus the inequality \eqref{sufficientimplicit} translates into
$$ \left[ \Tr \left(\frac{\partial W}{\partial t} \, \frac{\partial
W}{\partial y}\right) \right]^2 - 4 \det \left(\frac{\partial
W}{\partial t} \, \frac{\partial W}{\partial y}\right) \ge 0 \qquad
\text{at $(t_o,y_o)$.}$$ Now the equation $\Tr{A}^2 - 4 \det{A} \ge
0$ for a $2\times 2$ matrix $A$ is equivalent to the requirement
that $A$ has real eigenvalues. Thus we will have
\eqref{sufficientimplicit} if
\begin{equation}\label{derivativeproduct}
\frac{\partial W}{\partial t}(t_o,y_o) \, \frac{\partial W}{\partial
y}(t_o,y_o) =  J \Theta(t_o,y_o) \int_0^t W(\tau,y_o)^{-1} J \,
\frac{\partial \Theta}{\partial y}(\tau, y_o) W(\tau,y_o) \,
d\tau\end{equation} has real eigenvalues. To establish this, we
eliminate the $J$ matrices using the easily-proved formula $J Z =
(\det{Z}) (Z^T)^{-1} J$ along with $J^2 = -\id$ to obtain
$$\frac{\partial W}{\partial t}(t_o,y_o) \, \frac{\partial
W}{\partial y}(t_o,y_o) =\det{(\Theta(t_o,y_o))}
\Theta(t_o,y_o)^{-1} \int_0^{t_o} W(\tau,y_o)^T \frac{\partial
\Theta}{\partial y}(t_o,y_o) W(\tau,y_o) \, d\tau.$$

We now have the desired matrix as a product of a symmetric
positive-definite matrix and a symmetric matrix. It is easy to see
that such a product must have real eigenvalues: let $C$ be a
symmetric positive-definite matrix and $D$ a symmetric matrix. In an
eigenvector-basis of $C$, we have $C = \left(\begin{smallmatrix}
\lambda & 0 \\ 0 & \mu\end{smallmatrix}\right)$ and $D = \left(
\begin{smallmatrix} a & b \\ b& c\end{smallmatrix}\right)$, so that
$CD = \left(\begin{smallmatrix} \lambda a & \lambda b \\ \mu b & \mu
c\end{smallmatrix}\right)$. The expression $[\Tr{(CD)}]^2 - 4
\det{(CD)}$ becomes $(\lambda a - \mu c)^2 + 4\lambda \mu b^2$, and
this is nonnegative since $\lambda \mu > 0$. So $CD$ must have real
eigenvalues, and this establishes \eqref{sufficientimplicit}.
\end{proof}

\begin{corollary}\label{iwishthiswerereallyjustatheorem}
For each point $x_o\in M$, there is a family of $C^1$ functions
$f(\xi_o)$ defined on open subsets of the $2$-sphere, such that for
every $\xi_o$ there is a weakly epiconjugate point to the identity
located at $t=f(\xi_o)$. As a result, the set of all weakly
epiconjugate point locations obtainable by the technique of Theorem
\ref{main} consists of the closure of the union of open intervals.
\end{corollary}

\begin{proof}
For every component of the domain of $f$, the image of $f$. So we
have a family of intervals (possibly open, closed, or half-open) as
$\xi_o$ ranges over $S^2$. The actual set of epiconjugate points
obtained this way must be closed, since if we have an sequence of
conjugate point locations then any limit point must be an
epiconjugate point location. Hence we obtain some union of closed
intervals.
\end{proof}

By Theorem \ref{weakproposition} and Theorem
\ref{weakimpliesstrong}, we can conclude the following about the
intervals obtained from Corollary
\ref{iwishthiswerereallyjustatheorem}.
\begin{corollary}
All of the points in the nontrivial intervals obtained by Corollary
\ref{iwishthiswerereallyjustatheorem} are strongly epiconjugate.
\end{corollary}

\begin{proof}
Take any nontrivial closed interval of weakly epiconjugate point
locations $[a,b] \subset \mathbb{R}^+$. By Theorem
\ref{weakproposition}, there is a countable dense subset of this
interval consisting of weakly monoconjugate points. All other points
in the interval must therefore be strictly weakly epiconjugate, and
thus strongly epiconjugate. Finally, since a limit of strongly
epiconjugate points is also strongly epiconjugate, the entire
interval must consist of strongly epiconjugate points.
\end{proof}

It is not clear whether this technique actually generates more than
one interval; in the simple examples one can work out explicitly,
there is always just a single interval extending to infinity. The
following is a simple yet seemingly typical example.

\begin{proposition}
Suppose $u$ is a Killing field on $M$. (Such a $u$ is a steady
solution of the Euler equation, by Misio{\l}ek~\cite{m1}.) The
epiconjugate points along the corresponding geodesic $\eta$ form an
interval $[\tau,\infty)$ for some $\tau>0$.
\end{proposition}

\begin{proof}
Since $u$ is a Killing field, the metric pullback $\Lambda(t,x)$
must be the identity for all $t$ and $x$. Thus the equation
\eqref{itsnotreallyselfadjoint} takes the form
$$\frac{d^2\gamma}{dt^2} + \langle \omega_o(x_o),
\xi_o\rangle \xi_o \times \frac{d\gamma}{dt} = 0$$ with $\gamma(t)$
orthogonal to $\xi_o$ for all $t$. It is easy to see that this
equation has solutions vanishing at times
$$ \tau(x_o,\xi_o) = \frac{2\pi}{\langle \omega_o(x_o), \xi_o\rangle}$$
for any unit vector $\xi_o \in T_{x_o}M$. The infimum
$\tau(x_o)=\inf_{\xi_o\in S^2} \tau(x_o,\xi_o)$ occurs when $\xi_o$
is parallel to $\omega_o(x_o)$ and is $\tau(x_o) =
\frac{2\pi}{\lvert \omega_o(x_o)\rvert}$. The absolute minimum over
all $x_o \in M$ is
$$ \tau = \inf_{x_o\in M} \tau(x_o) = \frac{2\pi}{\lVert \omega_o\rVert_{L^{\infty}}}.$$
On the other hand, as $\xi_o$ approaches a vector orthogonal to
$\omega_o(x_o)$, the time $\tau(x_o,\xi_o)$ obviously approaches
infinity for any $x_o$.

Now it is also not hard to see that the actual first conjugate point
occurs at $\tau = \frac{2\pi}{\lVert \omega_o\rVert_{L^{\infty}}}$
(this is proved in \cite{p}). As a result, in this case Theorem
\ref{main} predicts all conjugate point locations.
\end{proof}

The behavior seen above for Killing fields appears to be typical;
for the flows whose conjugate points one can work out explicitly,
epiconjugate points occur in intervals extending to infinity. It is
conceivable that for nonsteady flows, the metric pullback $\Lambda$
may increase fast enough with time to prevent solutions of
\eqref{itsnotreallyselfadjoint} from vanishing at large times.
Obviously, if the solution of the Euler equation blows up at a
finite time, then all bets are off.

\bibliographystyle{plain}
\bibliography{Conjugate_Interval}

\end{document}